\documentclass[a4paper,11pt]{article}

\usepackage{latexsym,amssymb}
\usepackage{graphicx}
\usepackage{hyperref}
\usepackage{amsmath}
\usepackage[english]{babel}
\usepackage[matrix,arrow]{xy}

\usepackage{mathrsfs}

\def \be{\begin{eqnarray*}}
\def \ee{\end{eqnarray*}}
\def \ben{\begin{enumerate}}
\def \een{\end{enumerate}}
\def \beit{\begin{itemize}}
\def \eeit{\end{itemize}}
\def \bui#1#2{\mathrel{\mathop{\kern 0pt#1}\limits^{#2}}}
\def \buil#1#2{\mathrel{\mathop{\kern 0pt#1}\limits_{#2}}}
\def \bfll{\begin{flushleft}}
\def \efll{\end{flushleft}}
\def \bflr{\begin{flushright}}
\def \eflr{\end{flushright}}

\def \findemo{\hfill$\square$\\}
\def \fintdemo{\hfill$\surd$\\}
\def \lra{\longrightarrow}
\def \lmt{\longmapsto}
\def \ovl{\overline}

\def \wih{\widehat}
\def \wit{\widetilde}
\def \wnabla{\wit{\nabla}}
\def \cdotM{\buil{\cdot}{M}}
\def \cdotSi{\buil{\cdot}{\Sigma}}

\def \C{\mathbb{C}}

\def \R{\mathbb{R}}

\newcommand{\pa}[1]{\left(#1\right)}

\newtheorem{ethm}{Theorem}[section]

\newtheorem{erem}[ethm]{Remark}
\newtheorem{erems}[ethm]{Remarks}
\newtheorem{ecor}[ethm]{Corollary}
\newtheorem{prop}[ethm]{Proposition}

\def \C{\mathbb{C}}

\def \R{\mathbb{R}}

\title{A new upper bound for the Dirac operator on hypersurfaces}
\author{Nicolas Ginoux\footnote{Fakult\"at f\"ur Mathematik,
Universit\"at Regensburg,
D-93040 Regensburg,
E-mail: \texttt{nicolas.ginoux@mathematik.uni-r.de}},\, Georges Habib\footnote{Lebanese University, Faculty of Sciences II, Department of Mathematics, P.O. Box 90656 Fanar-Matn, Lebanon, 
E-mail: \texttt{ghabib@ul.edu.lb}} and Simon Raulot\footnote{Laboratoire de Math\'ematiques R. Salem
UMR $6085$ CNRS-Universit\'e de Rouen
Avenue de l'Universit\'e, BP.$12$
Technop\^ole du Madrillet
$76801$ Saint-\'Etienne-du-Rouvray, France. E-mail: \texttt{simon.raulot@univ-rouen.fr}}}

\begin{document}
\maketitle
\begin{minipage}{14cm}
{\it Dedicated to Oussama Hijazi for his sixtieth birthday and to Sebasti\'an Montiel}
\end{minipage}
\vspace{0.4cm}

\noindent\begin{center}\begin{tabular}{p{120mm}}
\begin{small}{\bf Abstract:} 
We prove a new upper bound for the first eigenvalue of the Dirac operator of a compact hypersurface in any Riemannian spin manifold carrying a non-trivial twistor spinor without zeros on the hypersurface.
The upper bound is expressed as the first eigenvalue of a drifting Schr\"odinger operator on the hypersurface.
Moreover, using a recent approach developed by O.~Hijazi and S.~Montiel, we completely characterize the equality case when the ambient manifold is the standard hyperbolic space.
\end{small}\\
\end{tabular}\end{center}


\section{Introduction}


Let $M^n\bui{\hookrightarrow}{\iota}\wit{M}^{n+1}$ be an oriented, compact (without boundary) and connected hypersurface of an $(n+1)$-dimensional Riemannian manifold $(\wit{M}^{n+1},g)$ equipped with the induced Riemannian metric also denoted by $g$.  

It is by now a well-known approach to use the min-max characterization of eigenvalues to derive upper bounds for the spectrum of differential operators on $M$ in terms of extrinsic geometric data. For example, if we consider the first positive eigenvalue $\lambda_1(\Delta)$ of the Laplace operator $\Delta:=-{\rm tr}_g({\rm Hess}_g)$ where ${\rm Hess}_g$ denotes the Hessian of $M$, a famous result of R.C.~Reilly \cite{Reilly77} states that if $\wit{M}$ is the Euclidean space $\mathbb{R}^{n+1}$ then 
\begin{eqnarray}\label{ReillyEuc}
\lambda_1(\Delta)\leq \frac{n}{{\rm Vol}(M)}\int_M H^2dv_g
\end{eqnarray}
where $H$ denotes the normalized mean curvature of $M$. The proof of this result uses in an essential way the Rayleigh characterization of $\lambda_1(\Delta)$ by choosing a modification of the coordinates functions as test functions. Moreover, it is a straightforward observation to see that equality occurs if and only if $M$ is a totally umbilical round sphere. As observed in \cite{ElSoufIlias}, this method directly applies for hypersurfaces in the unit sphere $\mathbb{S}^{n+1}$ leading to the counterpart of (\ref{ReillyEuc}) in this situation
\begin{eqnarray}\label{ReillySph}
\lambda_1(\Delta)\leq \frac{n}{{\rm Vol}(M)}\int_M\big(H^2+1\big)dv_g.
\end{eqnarray}
If the ambient manifold $\wit{M}$ is the standard hyperbolic space, there is also an optimal upper bound proved by A.~El Soufi and S.~Ilias \cite[Thm. 1]{ElSoufIlias} which improves a previous result of E.~Heintze \cite{Hein88} and which states that
\begin{eqnarray}\label{heintze}
\lambda_1(\Delta)\leq \frac{n}{{\rm Vol}(M)}\int_M\big(H^2-1)dv_g
\end{eqnarray}
with equality if and only if $M$ is a totally umbilical round sphere. 
All three estimates above follow actually from a much more general one, valid for submanifolds of any codimension and also proved by A.~El Soufi and S.~Ilias in \cite{ElSoufIlias}, assuming solely that the ambient manifold is conformally equivalent to an open subset of the sphere of the same dimension: under that assumption, they prove \cite[Thm. 2]{ElSoufIlias}
\begin{equation}\label{eq:ElSoufiIlias}
 \lambda_1(\Delta)\leq \frac{n}{{\rm Vol}(M)}\int_M\big(H^2+R(\iota))dv_g,
\end{equation}
where $R(\iota)$ is the normalized trace of the ambient sectional curvature on the tangent planes, see precise definition below (\ref{DefLPsi}).\\

Now if we assume the existence of a spin structure on $\wit{M}$ (which is the case for most classical ambient spaces), it induces a spin structure on the hypersurface $M$ and so we can define the spinor bundle $\Sigma M$ over $M$ as well as the associated Dirac operator $D_M$ (see Section \ref{s:prelim} and the references therein). When the ambient space $\wit{M}$ is the space form of constant sectional curvature $\kappa\in\{0,1,-1\}$, C.~B\"ar proved in \cite{Baer98} that 
\begin{eqnarray}\label{baer01}
\lambda_1(D_M^2)\leq\frac{n^2}{4{\rm Vol}(M)}\int_M\big(H^2+\kappa\big)dv_g
\end{eqnarray} 
if $\kappa=0,1$ and 
\begin{eqnarray}\label{baer-1}
\lambda_1(D_M^2)\leq\frac{n^2}{4}\sup_M\big(H^2+1\big)
\end{eqnarray}
for $\kappa=-1$. Here $\lambda_1(D_M^2)$ denotes the first non-negative eigenvalue of the square of the Dirac operator $D_M$ of $(M,g)$. Those estimates are consequences of the min-max characterization of $\lambda_1(D_M^2)$ and the fact that the space forms $\mathbb{R}^{n+1}$, $\mathbb{S}^{n+1}$ and $\mathbb{H}^{n+1}$ carry respectively parallel spinors, real and imaginary Killing spinors. In fact, taking the restriction of such a spinor field to the hypersurface as a test section in the Rayleigh quotient of $\lambda_1(D_M^2)$ gives immediately the previous inequalities. Note that these upper bounds hold for more general ambient manifolds since the proof only relies on the existence of one of such particular fields. For example, Inequality (\ref{baer01}) with $\kappa=0$ holds for compact oriented hypersurfaces in Calabi-Yau manifolds, hyper-K\"ahler and some other $7-$ and $8-$dimensional special Riemaniann manifolds. It also appears that both inequalities in (\ref{baer01}) are sharp since round geodesic spheres in the Euclidean space $\mathbb{R}^{n+1}$ and in the round sphere $\mathbb{S}^{n+1}$ satisfy the equality case. If $\kappa=0$, it has been recently proved by O.~Hijazi and S.~Montiel \cite{HijaziMontielhypsphere} that those are the only hypersurfaces for which equality is achieved.
The limiting case for hypersurfaces in the sphere seems to be out of reach at this time and could be considered as a spinorial analogue of the Yau conjecture about the first eigenvalue of the Laplace operator of minimal hypersurfaces in the unit sphere.
However, let us mention that there are non-minimal hypersurfaces in the sphere that satisfy the limiting case in (\ref{baer01}), see e.g. \cite{GinEgBaer03,GinSU2Q8}.\\

Regarding the proof of Inequality (\ref{baer-1}), it is not difficult to observe that there are no hypersurfaces which satisfy the equality case.
Modifying the computation of the Rayleigh quotient for $\lambda_1(D_M^2)$, the first named author improved this estimate into (see \cite[Thm. 1]{GinouxCRAS})
\begin{eqnarray}\label{ginoux-1}
\lambda_1(D_M^2)\leq\frac{n^2}{4}\sup_M\big(H^2-1\big),
\end{eqnarray}  
where equality occurs for totally umbilical round spheres in $\mathbb{H}^{n+1}$.
As we will see (Corollary \ref{t:eqcasehypspace}), those are in fact the only hypersurfaces for which Inequality (\ref{ginoux-1}) is an equality. 

\vspace{0.5cm}

In this paper, we prove a new upper bound for the first eigenvalue of the Dirac operator of $M$ when the ambient manifold $\wit{M}$ carries a twistor-spinor (Theorem \ref{t:lambda1L}).
This bound coincides with the first eigenvalue of an elliptic differential operator of order two whose definition depends among others on the norm of the twistor spinor along the hypersurface (see (\ref{DefLPsi})) and which belongs to a particular class of operators: the drifting Schr\"odinger operators, that is, of the form drifting Laplacian plus potential (see Remark \ref{r:driftinglaplacian}). It is important to note that this estimate contains all the (up to date) known upper estimates {\it \`a la Reilly} (see Remark \ref{r:knownestimates}).  In a second part, we adapt the approach developed by O.~Hijazi and S.~Montiel \cite{HijaziMontielhypsphere} to prove that, assuming the existence of imaginary Killing spinors for two opposite constants on $\wit{M}$, the only hypersurfaces satisfying the equality case in our previous estimate are the totally umbilical ones (Theorem \ref{t:eqcasegeneral}). 
In particular, only the geodesic hyperspheres satisfy that limiting case in the hyperbolic space (Corollary \ref{t:eqcasehypspace}).
We also examine the setting of pseudo-hyperbolic spaces (see Corollary \ref{c:Baumclass}). 
 

\section{Preliminaries and notations}\label{s:prelim}


In this section, we briefly introduce the geometric setting and fix the notations of this paper. For more details on those preliminaries we refer for example to \cite{LM89}, \cite{Friedrichlivre} or \cite[Ch. 1]{Ginouxbook}. 

We consider $M^n\bui{\hookrightarrow}{\iota}\wit{M}^{n+1}$ an oriented $n$-di\-men\-sio\-nal Riemannian hypersurface with $n\geq 2$, isometrically immersed into an $(n+1)$-dimensional Riemannian spin manifold $(\wit{M}^{n+1},g)$ with a fixed spin structure. We denote by $\nu$ the unit inner normal vector field induced by both orientations, that is, such that $(E_1,\cdots,E_n,\nu_x)$ is an oriented basis of $T_x\wit{M}_{|_M}$ if and only if $(E_1,\cdots,E_n)$ is an oriented basis of $T_xM$ for $x\in M$. We endow $M$ with the spin structure induced by the one on $\wit{M}$ and let $\Sigma M\to M$ denotes the associated spinor bundle. Setting 
\[\Sigma:=\left\{\begin{array}{ll}\Sigma M&\textrm{if }n\textrm{ is even}\\\Sigma M\oplus\Sigma M&\textrm{if }n\textrm{ is odd,}\end{array}\right.\]
the bundles $\Sigma$ and the restriction $\Sigma\wit{M}_{|_M}$ to $M$ of the spinor bundle of $\wit{M}$ can be identified such that
\beit
\item both natural Hermitian inner products -- that we hence denote by $\langle\cdot\,,\cdot\rangle$ -- coincide,

\item the Clifford multiplication ``\,$\cdot$\,'' on $\wit{M}$ and ``\,$\cdotM$\,'' on $M$ are related by 
\begin{equation}\label{d:CliffMult}
X\cdotSi\,:=X\cdot\nu\cdot\simeq
\left\{\begin{array}{ll}X\cdotM&\textrm{if }n\textrm{ is even}\\ 
X\cdotM\oplus-X\cdotM&\textrm{if }n\textrm{ is odd,}
\end{array}\right.
\end{equation}
for all $X\in TM$,

\item the spin Levi-Civita connections  $\wnabla$ on $\Sigma\wit{M}$ and $\nabla$ on $\Sigma$ are related by the spin Gau{\ss} formula 
\begin{equation}\label{eq:Gaussspin}
\wnabla_X\varphi=\nabla_X\varphi+\frac{A(X)}{2}\cdot\nu\cdot\varphi, 
\end{equation}
for all $X\in\Gamma(TM)$ and $\varphi\in\Gamma(\Sigma)$. Here $A:=-\wnabla\nu$ denotes the Weingarten map of the immersion.
\eeit
The extrinsic Dirac operator of $M$ is the first order elliptic differential operator of order one acting on sections of $\Sigma$ locally given by
\begin{eqnarray*}
D:=\sum_{j=1}^ne_j\cdot\nu\cdot\nabla_{e_j}.
\end{eqnarray*}
It is a well-known fact that it defines an essentially self-adjoint operator with respect to the $L^2$-scalar product on $\Sigma$ so that if $M$ is compact, its spectrum is an unbounded sequence of real numbers.
By convention and in the whole article, the spectrum ${\rm spec}(P)$ with multiplicities of a given elliptic self-adjoint operator $P$ will be denoted by a sequence $\left(\lambda_k(P)\right)_{k\geq1}$, with the convention that $\lambda_1(P)$ is the smallest eigenvalue if ${\rm spec}(P)$ is bounded below and is the smallest \emph{nonnegative} eigenvalue otherwise.\\

With respect to the previous identifications, the Dirac operator $D$ is nothing but the Dirac operator $D_M$ of $(M,g)$ if $n$ is even and $D_M\oplus -D_M$ if $n$ is odd, so that studying the spectrum of the intrinsic Dirac operator $D_M$ for the spin Riemannian structure induced on the hypersurface $M$ is equivalent to study the spectrum of the extrinsic Dirac operator $D$ on the hypersurface $M$. It is also relevant here to recall that the commutator of $D$ and $D^2$ with functions are given by
\begin{equation}\label{eq:commutDf}
D(f\varphi)=fD\varphi+\nabla f\cdot\nu\cdot\varphi 
\end{equation}
and
\begin{equation}\label{eq:commutD2f}
D^2(f\varphi)=fD^2\varphi-2\nabla_{\nabla f}\varphi+(\Delta f)\varphi,
\end{equation}
for all $f\in C^\infty(M)$ and $\varphi\in\Gamma(\Sigma)$. Here $H:=(1/n){\rm tr}(A)$ denotes the mean curvature function of $M$ in $\wit{M}$.

Another operator of particular interest in this work is the Dirac-Witten operator $\wih{D}$ on $M$. It is also a first order elliptic operator acting on the restricted spinor bundle $\Sigma$ and locally defined by $\wih{D}:=\sum_{j=1}^ne_j\cdot\wnabla_{e_j}$. It is related to the extrinsic Dirac operator by the following formula
\begin{equation}\label{eq:GaussDirac}
D\varphi=-\nu\cdot\wih{D}\varphi+\frac{nH}{2}\varphi 
\end{equation}
and to its squared by
\begin{equation}\label{eq:GaussDirac2}
D^2\varphi=\wih{D}^2\varphi+\frac{n^2H^2}{4}\varphi+\frac{n}{2}\nabla H\cdot\nu\cdot\varphi, 
\end{equation}
for every $\varphi\in\Gamma(\Sigma)$.


\section{Upper bounds in terms of a Laplace-type operator}\label{s:upperblambda1L}


In this section, we prove a new upper bound for the smallest eigenvalue of the squared Dirac operator $D^2$ when the ambient manifold $\wit{M}$ is endowed with a twistor spinor. Recall that a \emph{twistor spinor} on a Riemannian spin manifold $(\wit{M}^{n+1},g)$ is a section $\psi\in\Gamma(\Sigma\wit{M})$ satisfying
\begin{eqnarray}\label{d:twistspin}
\wnabla_X\psi=-\frac{1}{n+1}X\buil{\cdot}{\wit{M}}D_{\wit{M}}\psi
\end{eqnarray}
for all $X\in\Gamma(T\wit{M})$. Here $D_{\wit{M}}$ represents the Dirac operator of $\wit{M}$.
Non-zero twistor-spinors have a discrete vanishing set and only exist for particular conformal classes (see for example the standard reference \cite{BFGK} or \cite[App. A]{Ginouxbook} for a short account). It should also be pointed out that parallel spinors, real and imaginary Killing spinors are twistor spinors which are, in addition, eigensections for the Dirac operator $D_{\wit{M}}$ respectively associated to the eigenvalue zero, or to real or purely imaginary eigenvalues. They exist on each simply connected complete space form of constant curvature. Assume now that such a spinor field $\psi$ is given on $\wit{M}$ and also assume that it has no zero on the hypersurface $M$. We define the differential operator $L_\psi$ acting on smooth functions on $M$ by
\begin{eqnarray}\label{DefLPsi}
L_\psi f:=\Delta f-2g(\nabla\ln|\psi|,\nabla f)+\frac{n^2}{4}(H^2+R(\iota))f.
\end{eqnarray}
for $f\in C^\infty(M)$. Here $R(\iota):=\frac{1}{n(n-1)}\left(\wit{S}-2\,\wit{\rm ric}(\nu,\nu)\right)$, $\wit{S}$ and $\wit{\rm ric}$ are respectively the scalar curvature and the Ricci tensor (seen as a symmetric $2$-tensor) of the manifold $\wit{M}$. Although this operator is not symmetric with respect to the $L^2$-scalar product on $(M^n,g)$, we observe that it has the following interesting analytic properties:
\begin{prop}\label{AnalyticL}
The operator $L_\psi$ is elliptic and if $M$ is closed, it is self-adjoint with respect to the $L^2$-scalar product on $(M^{n},\ovl{g}:=|\psi|^{\frac{4}{n}}\,g)$.
\end{prop}

{\it Proof:}
Since $L_\psi$ is of second order and its leading part is the scalar Laplacian, it is clearly elliptic.
Because of $\ovl{g}=|\psi|^{\frac{4}{n}}\,g$, we have $dv_{\ovl{g}}=|\psi|^2\,dv_g$ and we can write for any $f,h\in C^\infty(M)$:
\[\int_M(L_\psi f)h\,dv_{\ovl{g}}=\int_M\left(\Delta f-2g(\nabla\ln|\psi|,\nabla f)+\frac{n^2}{4}\big(H^2+R(\iota)\big)f\right)h|\psi|^2dv_g.\]
Performing a partial integration, we have for the first term
\be 
\int_M(\Delta f)h|\psi|^2dv_g&=&\int_Mg(\nabla f,\nabla h)|\psi|^2+g(\nabla f,\nabla(|\psi|^2))hdv_g\\
&=&\int_Mg(\nabla f,\nabla h)|\psi|^2+2g(\nabla f,\nabla\ln|\psi|)h|\psi|^2dv_g.
\ee
Therefore, the first-order term in $\nabla\ln|\psi|$ simplifies and we obtain
\[\int_M(L_\psi f)h\,dv_{\ovl{g}}=\int_M\left(g(\nabla f,\nabla h)+\frac{n^2}{4}\big(H^2+R(\iota)\big)fh\right)|\psi|^2dv_g,\]
which is clearly symmetric in $(f,h)$.
This implies that $L_\psi$ is formally self-adjoint with respect to the metric $\ovl{g}$.
Since $M$ is closed, we conclude that $L_\psi$ is essentially self-adjoint in $L^2(M)$.
\findemo

\begin{erems}\label{r:driftinglaplacian}
\noindent\begin{enumerate}
\item {\rm The operator $L_\psi$ defined in (\ref{DefLPsi}) is of the form \emph{drifting Laplacian} (also called \emph{Laplacian with drift}, \emph{Bakry-Emery Laplacian}, \emph{weighted Laplacian} or \emph{Witten Laplacian} in the literature) plus potential, this is the reason we refer to these operators as \emph{drifting Schr\"odinger operators}. Indeed, a drifting Laplacian is an operator of the form 
$$
\begin{array}{cll}
C^\infty(M) & \bui{\lra}{L_h} & C^\infty(M)\\
f & \mapsto & \Delta f-g(\nabla h,\nabla f)
\end{array}
$$
for some function $h\in C^\infty(M)$. It is elliptic and self-adjoint with respect to the measure $e^hd\mu_g$. Actually, a drifting Laplacian is always \emph{unitarily equivalent} to a Schr\"odinger operator: in the notations above, the operator $L_h$ is unitarily equivalent to $\Delta-\frac{1}{2}\Delta h+\frac{1}{4}|\nabla h|_g^2$ (see for example \cite[p.28]{Setti98}).
} 
\item {\rm Note that if $|\psi|$ is constant on $M$ (which is the case if $\psi$ is either a parallel or a real Killing spinor on $\wit{M}$), then the operator 
\begin{eqnarray*}
L_\psi=\Delta+\frac{n^2}{4}(H^2+R(\iota)) 
\end{eqnarray*}
does not depend on $\psi$.}
\end{enumerate}
\end{erems}

Proposition \ref{AnalyticL} implies that the spectrum of $L_\psi$ is purely discrete.
We will denote by $\lambda_1(L_\psi)$ its first eigenvalue, which satisfies the following min-max characterization
\begin{eqnarray}\label{RayleL}
\lambda_1(L_\psi)=\underset{f\in C^\infty(M)\setminus\{0\}}{\inf}\left(\frac{\int_Mf(L_\psi f)dv_{\ovl{g}}}{\int_Mf^2dv_{\ovl{g}}}\right).
\end{eqnarray}
 
We are now ready to give the precise statement of the first main result of this paper, namely:
\begin{ethm}\label{t:lambda1L}
Assume $M$ is a closed oriented hypersurface isometrically immersed in a Riemannian spin manifold $(\wit{M}^{n+1},g)$. If there exists a non trivial twistor-spinor $\psi$ on $\wit{M}$ with $\psi_x\neq0$ for all $x\in M$ then we have
\begin{equation}\label{eq:majlambda1L}
\lambda_1(D_M^2)\leq\lambda_1(L_\psi). 
\end{equation}
\end{ethm}

{\it Proof}:
We apply the min-max characterization of $\lambda_1(D_M^2)=\lambda_1(D^2)$ using $f\psi$ as a test-section, where $L_\psi f=\lambda_1(L_\psi)f$. The following computations rely on a large extent on those in the proof of \cite[Thm. 5.2.3]{Ginouxbook}. 

First, if $f\in C^\infty(M)$ is an arbitrary smooth function on $M$, then using (\ref{eq:commutD2f}), (\ref{eq:GaussDirac2}), (\ref{eq:Gaussspin}) and the fact that $\psi$ is a twistor-spinor on $\wit{M}$, one obtains 
\begin{eqnarray}\label{eq:D2fpsisptw1}
\nonumber D^2(f\psi)&\bui{=}{\rm(\ref{eq:commutD2f})}&fD^2\psi-2\nabla_{\!\nabla f}\psi+(\Delta f)\psi\\
\nonumber &\bui{=}{\rm(\ref{eq:GaussDirac2})}&f\big(\wih{D}^2\psi+\frac{n^2H^2}{4}\psi+\frac{n}{2}\nabla H\cdot\nu\cdot\psi\big)-2\nabla_{\!\nabla f}\psi+(\Delta f)\psi\\
\nonumber &\bui{=}{\rm(\ref{eq:Gaussspin})}&f(\wih{D}^2\psi+\frac{n^2H^2}{4}\psi+\frac{n}{2}\nabla H\cdot\nu\cdot\psi)\\
\nonumber & &-2(\wnabla_{\!\nabla f}\psi-\frac{A(\nabla f)}{2}\cdot\nu\cdot\psi)+(\Delta f)\psi\\
\nonumber &=&f\big(\wih{D}^2\psi+\frac{n^2H^2}{4}\psi+\frac{n}{2}\nabla H\cdot\nu\cdot\psi\big)+\frac{2}{n+1}\nabla f\cdot D_{\wit{M}}\psi\\
& &+A(\nabla f)\cdot\nu\cdot\psi+(\Delta f)\psi.
\end{eqnarray}
Next we compute $\wih{D}^2\psi$, using again the fact that $\psi$ is a twistor-spinor, fact which implies in particular the following identity (see e.g. \cite[Prop. A.2.1]{Ginouxbook}):
\begin{equation}\label{eq:nablaDpsisptw}
\wnabla_{\! X}(D_{\wit{M}}\psi)=\frac{n+1}{n-1}\Big(-\frac{1}{2}\wit{\mathrm{Ric}}(X)\cdot\psi+\frac{\wit{S}}{4n}X\cdot\psi\Big),
\end{equation}
for every $X\in\Gamma(T\wit{M})$ and where $\wit{\mathrm{Ric}}$ denotes the Ricci tensor of $(\wit{M}^{n+1},g)$ (seen as an endomorphism of the tangent bundle of $\wit{M}$). Thus we have
\begin{eqnarray}\label{eq:wihD2psisptw}
\nonumber\wih{D}^2\psi&=&\wih{D}(\sum_{j=1}^n e_j\cdot\wnabla_{e_j}\psi)\\
\nonumber&\bui{=}{(\rm\ref{d:twistspin})}& \frac{n}{n+1}\wih{D}(D_{\wit{M}}\psi)\\
\nonumber&\bui{=}{(\rm\ref{eq:nablaDpsisptw})}&\frac{n}{n-1}\sum_{j=1}^n\big(-\frac{1}{2}e_j\cdot\wit{\mathrm{Ric}}(e_j)\cdot\psi+\frac{\wit{S}}{4n}e_j\cdot e_j\cdot\psi\big)\\
\nonumber&=&\frac{n}{n-1}\Big(\frac{\wit{S}}{2}\psi+\frac{1}{2}\nu\cdot\wit{\mathrm{Ric}}(\nu)\cdot\psi-\frac{\wit{S}}{4}\psi\Big)\\
\nonumber&=&\frac{n}{n-1}\Big(\frac{n(n-1)}{4}R(\iota)\psi+\frac{1}{2}\nu\cdot\wit{\mathrm{Ric}}(\nu)^{\mathrm{T}}\cdot\psi\Big)\\
&=&\frac{n^2}{4}R(\iota)\psi+\frac{n}{2(n-1)}\nu\cdot\wit{\mathrm{Ric}}(\nu)^{\mathrm{T}}\cdot\psi,
\end{eqnarray}
where $\wit{\mathrm{Ric}}(\nu)^{\mathrm{T}}:=\sum_{j=1}^n\wit{\mathrm{ric}}(\nu,e_j)e_j$ denotes the tangential projection of $\wit{\mathrm{Ric}}(\nu)$ on $TM$. Combining (\ref{eq:D2fpsisptw1}) with (\ref{eq:wihD2psisptw}), we deduce that 
\begin{eqnarray}\label{eq:D2fpsisptw2}
\nonumber D^2(f\psi)&=&\frac{n^2}{4}\big(H^2+R(\iota)\big)f\psi+\frac{nf}{2}\nabla H\cdot\nu\cdot\psi+\frac{nf}{2(n-1)}\nu\cdot\wit{\mathrm{Ric}}(\nu)^{\mathrm{T}}\cdot\psi\\
& &+\frac{2}{n+1}\nabla f\cdot D_{\wit{M}}\psi+A\big(\nabla f\big)\cdot\nu\cdot\psi+(\Delta f)\psi.
\end{eqnarray}
Using again that $\psi$ is a twistor-spinor on $(\wit{M}^{n+1},g)$, we obtain that for every $f\in C^\infty(M)$,
\be\Re e\big(\langle D^2(f\psi),f\psi\rangle\big)&\bui{=}{(\rm\ref{eq:D2fpsisptw2})}&\frac{n^2}{4}\big(H^2+R(\iota)\big)f^2|\psi|^2+\frac{2f}{n+1}\Re e\big(\langle \nabla f\cdot D_{\wit{M}}\psi,\psi\rangle\big)\\
& &+f(\Delta f)|\psi|^2\\
&=&\frac{n^2}{4}\big(H^2+R(\iota)\big)f^2|\psi|^2-g\big(f\nabla f,\nabla(|\psi|^2)\big)\\
& &+f(\Delta f)|\psi|^2\\
&=&f\,\Big(\Delta f-2g(\nabla f,\nabla\ln|\psi|)+\frac{n^2}{4}\big(H^2+R(\iota)\big)f\Big)|\psi|^2\\
&=&f(L_\psi f)|\psi|^2.
\ee
The min-max principle for $\lambda_1(D^2)$ implies that, for any $f\in C^\infty(M)\setminus\{0\}$,
\be 
\lambda_1(D^2)\leq\frac{\int_M\Re e\big(\langle D^2(f\psi),f\psi\rangle\big)dv_g}{\int_M|f\psi|^2dv_g}
=\frac{\int_Mf(L_\psi f)dv_{\ovl{g}}}{\int_Mf^2dv_{\ovl{g}}},
\ee
therefore, 
\[\lambda_1(D^2)\leq\buil{\inf}{f\in C^\infty(M,\R)\setminus\{0\}}\left(\frac{\int_Mf(L_\psi f)dv_{\ovl{g}}}{\int_Mf^2dv_{\ovl{g}}}\right)\]
which from (\ref{RayleL}) gives the inequality (\ref{eq:majlambda1L}). 
\findemo

\begin{erems}\label{r:knownestimates}
{\rm\noindent\ben\item The estimate (\ref{eq:majlambda1L}) contains all known upper estimates \emph{\`a la Reilly} for $\lambda_1(D_M^2)$.
Indeed, we observe that by taking $f=1$ in the Rayleigh quotient of $L_\psi$, we have 
\begin{eqnarray*}
\lambda_1(L_\psi)\leq\frac{n^2}{4{\rm Vol}(M)}\int_M\big(H^2+R(\iota)\big)dv_g
\end{eqnarray*}
if $|\psi|$ is constant and
\[\lambda_1(L_\psi)\leq\frac{n^2}{4}\sup_M\big(H^2+R(\iota)\big)\]
otherwise.
Those give exactly the inequalities (\ref{baer01}) by C.~B\"ar in \cite{Baer98} and (\ref{ginoux-1}) in \cite{GinouxCRAS} by the first named author. On the other hand, for $f=|\psi|^{-1}$ (w.r.t. the metric $\ovl{g}$ defined above) we deduce that
\[\lambda_1(L_\psi)\leq\frac{n^2}{4{\rm Vol}(M)}\int_M\big(H^2+R(\iota)\big)dv_g+\frac{1}{{\rm Vol}(M)}\int_M\big|d\ln|\psi|\big|^2dv_g\]
which was proved by the first-named author \cite[Thm. 1]{Gin02}.
\item It is interesting to compare (\ref{eq:majlambda1L}) with (\ref{eq:ElSoufiIlias}).
On the one hand, we do not obtain in the spinorial setting the exact analogue of (\ref{eq:ElSoufiIlias}) for $\wit{M}$ conformally equivalent to an open subset of the sphere $\mathbb{S}^{n+1}$. 
Of course, this must be expected since otherwise in dimension $2$ this would mean that the Willmore functional bounds $\lambda_1(D_M^2)\cdot{\rm Area}(M^2,g)$ from above; but there is no conformal upper bound for the smallest positive Dirac eigenvalue on unit-area-metrics, as shown in \cite[Thm. 1.1]{AmmJammes}.
Note that this does not prevent the analogue of (\ref{heintze}) to possibly hold true for the Dirac operator, which is still an open question.
On the other hand, our assumption on $\wit{M}$ in Theorem \ref{t:lambda1L} is much more general since not only open subsets of spheres with conformal metrics allow twistor-spinors.
We refer to \cite{KR98} for the classification of Riemannian spin manifolds with twistor-spinors. 
\een}
\end{erems}

We now look at the equality case of the previous estimate in the case of the twistor spinor is also an eigenspinor for the Dirac operator of $\wit{M}$. More precisely, we prove:
\begin{prop}\label{p:equality}
Under the same assumptions as in {\rm Theorem \ref{t:lambda1L}}, assume moreover that equality is achieved in {\rm (\ref{eq:majlambda1L})}.
Then
\ben
\item[$1.$] if $\psi$ is a parallel spinor on $\wit{M}^{n+1}$, one has 
\begin{eqnarray*}
A(\nabla\ln|f|)=-\frac{n}{2}\nabla H
\end{eqnarray*} 
for any eigenfunction $f$ of $L_\psi$ associated with $\lambda_1(L_\psi)$;

\item[$2.$] if $\psi$ is a real (resp. imaginary) Killing spinor on $\wit{M}=\mathbb{S}^{n+1}$ (resp. $\wit{M}=\mathbb{H}^{n+1}$), the mean curvature $H$ is constant and in particular $\lambda_1(D_M^2)=\frac{n^2}{4}\left(H^2+\kappa\right)$.
\een
\end{prop}

{\it Proof:}
\ben
\item If (\ref{eq:majlambda1L}) is an equality and $\psi$ is a parallel spinor, then the min-max principle yields $D^2(f\psi)=\lambda_1(D^2)f\psi$ for any eigenfunction $f$ of $L_\psi$ associated with $\lambda_1(L_\psi)=\lambda_1(D^2)$. But (\ref{eq:D2fpsisptw2}) together with $\wit{\rm Ric}=0$ and $D_{\wit{M}}\psi=0$ (both provided by $\wnabla\psi=0$) implies
\be
\lambda_1(D^2)f\psi&=&\frac{n^2H^2}{4}f\psi+\frac{nf}{2}\nabla H\cdot\nu\cdot\psi+A\big(\nabla f\big)\cdot\nu\cdot\psi+(\Delta f)\psi\\
&=&(L_\psi f)\psi+\big(A\big(\nabla f\big)+\frac{nf}{2}\nabla H \big)\cdot\nu\cdot\psi.
\ee
With $\lambda_1(D^2)=\lambda_1(L_\psi)$, we deduce that 
\begin{eqnarray*}
\Big(A\big(\nabla f\big)+\frac{nf}{2}\nabla H\Big)\cdot\nu\cdot\psi=0
\end{eqnarray*}
which, since $\psi\neq0$, gives $A\big(\nabla f\big)+\frac{nf}{2}\nabla H=0$. Since any eigenfunction for $L_\psi$ associated with the eigenvalue $\lambda_1(L_\psi)$ is either positive or negative, we easily conclude.

\item Assume first $\wit{M}^{n+1}$ carries real Killing spinors and let $\psi$ be a non-zero $(\varepsilon/2)$-Killing spinor for some $\varepsilon\in\{\pm1\}$, that is, $\wnabla_X\psi=(\varepsilon/2)X\cdot\psi$ for all $X\in\Gamma(T\wit{M})$. Again, one obtains $D^2(f\psi)=\lambda_1(D^2)f\psi$ for any eigenfunctions $f\in C^\infty(M)$ associated to $\lambda_1(L_\psi)$. Fixing such an $f$, the identity (\ref{eq:D2fpsisptw2}) yields
\be
\lambda_1(D^2)f\psi & = & (L_\psi f)\psi+\Big(A\big(\nabla f\big)+\frac{nf}{2}\nabla H\Big)\cdot\nu\cdot\psi-\varepsilon\nabla f\cdot\psi.
\ee 
With $\lambda_1(D^2)=\lambda_1(L_\psi)$, we deduce that
\begin{eqnarray*} 
\Big(A\big(\nabla f\big)+\frac{nf}{2}\nabla H\Big)\cdot\nu\cdot\psi-\varepsilon\nabla f\cdot\psi=0.
\end{eqnarray*}
In particular, denoting $Y_\varepsilon:=-\varepsilon\nabla f$ and $X:=A(\nabla f)+\frac{nf}{2}\nabla H$, we have $(Y_\varepsilon+X\wedge\nu)\cdot\psi=0$. At this point, we need the following claim:\\

{\bf Claim}: {\sl Let $\alpha\in\Lambda^*\R^{n+1}\otimes\C$. If $n$ is odd, then $\delta_{n+1}(\alpha)=0$ if and only if $\alpha=0$. If $n$ is even, then the same equivalence holds for $\alpha\in\Lambda^*\R^{n}\otimes\C$.}\\

{\it Proof of Claim}: 
Recall that the spinor representation $\delta_k:\mathbb{C}\mathrm{l}_k\lra\mathrm{End}_{\C}(\Sigma_k)$ of the complex Clifford algebra in dimension $k$ is a complex-linear isomorphism for $k$ even (but obviously not for $k$ odd). So if $n$ is odd, the claim follows directly from this fact. If $n$ is even and $\alpha\in\Lambda^*\R^{n}\otimes\C$, then $\Sigma_n\cong\Sigma_{n+1}$ and it is a simple trick to rewrite $\delta_{n+1}(\alpha)$ under the form $\delta_n(\check{\alpha})$ for a form $\check{\alpha}\in\Lambda^*\R^n\otimes\C$ having the same coefficients as $\alpha$ in the canonical basis of $\Lambda^*\R^n\otimes\C$ up to sign and some power of $i$. Namely, write 
\begin{eqnarray*}
\alpha=\buil{\sum}{1\leq j_1<\ldots<j_k\leq n}\alpha_{j_1,\ldots,j_k}e_{j_1}^*\wedge\ldots\wedge e_{j_k}^*,
\end{eqnarray*}
where $(e_1,\ldots,e_n,e_{n+1})$ is the canonical basis of $\R^{n+1}$.
Let $\omega_n^{\C}$ denotes the complex volume form on $\R^n$ as defined in the proof of Proposition \ref{p:equality}, which acts on $\Sigma_n$ via $\delta_n(\omega_n^{\C})=\mathrm{Id}_{\Sigma_n^+}\oplus-\mathrm{Id}_{\Sigma_n^-}$. Since $\delta_{n+1}(ie_{n+1})=\delta_n(\omega_n^{\C})$ and $\delta_n(v)=\delta_{n+1}(v)\circ\delta_{n+1}(e_{n+1})$ for all $v\in\R^n$, we have after some calculations
\be 
\delta_{n+1}(\alpha)&=&\sum_{\stackrel{1\leq j_1<\ldots<j_k\leq n}{k\textrm{ even}}}\alpha_{j_1,\ldots,j_k}\delta_{n}(e_{j_1})\circ\ldots\circ \delta_{n}(e_{j_k})\\
&&+i\sum_{\stackrel{1\leq j_1<\ldots<j_k\leq n}{k\textrm{ odd}}}\alpha_{j_1,\ldots,j_k}\delta_{n}(e_{j_1})\circ\ldots\circ \delta_{n}(e_{j_k})\circ \delta_n(\omega_n^{\C}).
\ee
Now it is an elementary computation to show that, for any $\beta\in\Lambda^k\R^n$, one has $\delta_n(\beta)\circ\delta_n(e_1^*\wedge\ldots\wedge e_n^*)=(-1)^{\frac{k(k+1)}{2}}\delta_n(*\beta)$, where $*:\Lambda^*\R^n\to\Lambda^*\R^n$ is the Hodge-star operator.
Therefore, we obtain
\be 
\delta_{n+1}(\alpha)&=&\sum_{\stackrel{1\leq j_1<\ldots<j_k\leq n}{k\textrm{ even}}}\alpha_{j_1,\ldots,j_k}\delta_{n}(e_{j_1})\circ\ldots\circ \delta_{n}(e_{j_k})\\ 
& & +c_{n,k}\sum_{\stackrel{1\leq j_1<\ldots<j_k\leq n}{k\textrm{ odd}}}\alpha_{j_1,\ldots,j_k}\delta_{n}(*(e_{j_1}^*\wedge\ldots\wedge e_{j_k}^*))\\
&=&\delta_n(\check{\alpha}),
\ee
where we let $c_{n,k}:=i^{\frac{n}{2}+1}(-1)^{\frac{k(k+1)}{2}}$ and
\[\check{\alpha}:=\sum_{\stackrel{1\leq j_1<\ldots<j_k\leq n}{k\textrm{ even}}}\alpha_{j_1,\ldots,j_k}e_{j_1}^*\wedge\ldots\wedge e_{j_k}^*+c_{n,k}\sum_{\stackrel{1\leq j_1<\ldots<j_k\leq n}{k\textrm{ odd}}}\alpha_{j_1,\ldots,j_k}*(e_{j_1}^*\wedge\ldots\wedge e_{j_k}^*).\]
As a consequence, if $\delta_{n+1}(\alpha)\sigma=0$ for all $\sigma\in\Sigma_{n+1}\cong\Sigma_n$, then $\delta_n(\check{\alpha})=0$ and the fact mentioned above implies $\check{\alpha}=0$; since $n$ is even, each form $*(e_{j_1}^*\wedge\ldots\wedge e_{j_k}^*)$ is of odd degree when $k$ is odd and therefore $\alpha_{j_1,\ldots,j_k}=0$ for all $1\leq j_1<\ldots<j_k\leq n$, that is, $\alpha=0$.
This concludes the proof of the claim.
\fintdemo

If $\wit{M}^{n+1}$ is isometric to the standard round sphere $\mathbb{S}^{n+1}$, then it carries a \emph{maximal} number (that is $2^{[\frac{n+1}{2}]}$) of linearly independent $(\varepsilon/2)$-Killing spinors, then $(Y_\varepsilon+X\wedge\nu)\cdot\psi=0$ holds pointwise for every $\psi\in\Sigma_x\wit{M}$. If $n$ is odd, then the claim yields $Y_\varepsilon+X\wedge\nu=0$, which implies $X=Y_\varepsilon=0$, that is, $f$ and $H$ are constant. If $n$ is even, one may rewrite 
\begin{eqnarray*}
Y_\varepsilon\cdot\psi+X\cdot\nu\cdot\psi=iY_\varepsilon\cdot i\nu\cdot \nu\cdot\psi+X\cdot\nu\cdot\psi=(X-iY_\varepsilon\lrcorner\omega_M^{\C})\cdot\nu\cdot\psi, 
\end{eqnarray*}
where $\omega_M^{\C}:=i^{[\frac{n+1}{2}]}e_1^*\wedge\ldots\wedge e_n^*\in\Gamma(\Lambda^nT^*M\otimes\C)$ is the complex volume form on $M$. Again, the claim yields $X-iY_\varepsilon\lrcorner\omega_M^{\C}=0$.
If $n>2$, then comparing the degrees yields $X=Y_\varepsilon=0$, that is, $f$ and $H$ are constant.
If $n=2$, then an elementary computation gives $Z\lrcorner\omega_M^{\C}=iJ(Z)$ for every $Z\in\Gamma(TM)$, where $J$ is the K\"ahler structure associated to the metric and the orientation on $(M^2,g)$.
In that case, one obtains $X+J(Y_\varepsilon)=0$. However on the standard sphere $\mathbb{S}^3$, both spaces of $\frac{1}{2}$- and $-\frac{1}{2}$-Killing spinors have maximal dimension $2$, therefore $X+J(Y_\varepsilon)=0$ for both $\varepsilon\in\{\pm1\}$, which implies $X=Y_\varepsilon=0$ and hence $f$ and $H$ are constant.\\

The case of imaginary Killing spinors is much the same up to replacing $\varepsilon$ by $i\varepsilon$.
One obtains at the end $(iY_\varepsilon+X\wedge\nu)\cdot\psi=0$ for all $(i\varepsilon/2)$-Killing spinors $\psi$ on $\wit{M}^{n+1}$. The same arguments as above lead to $X=Y_\varepsilon=0$. Remark that in the case $n=2$, one does not need the existence of maximal spaces of $\frac{i\varepsilon}{2}$-Killing spinors for both $\varepsilon\in\{\pm1\}$ since $X$ and $Y_\varepsilon$ are real vector fields on $M$.
\een
\findemo

\begin{erem}
{\rm It is quite surprising that in the case where $\psi$ is a parallel spinor we cannot conclude that the mean curvature of $M$ must be constant. In fact, we are left to prove that if there exists a smooth positive function $f\in C^\infty(M)$ such that
\begin{eqnarray*}
\Delta f+\frac{n^2H^2}{4}f=\lambda_1(D)^2f\quad{\rm and}\quad A(\nabla\ln f)=-\frac{n}{2}\nabla H
\end{eqnarray*}
then $f$ (or, equivalently, $H$) is constant on $M$.
}
\end{erem}



\section{Equality case in presence of imaginary Killing spinors}\label{s:eqcasehypspace}


In this section, we focus on the equality case of our estimate (\ref{eq:majlambda1L}) when the ambient manifold $\wit{M}$ carries an imaginary Killing spinor. According to Proposition \ref{p:equality}, it also corresponds to the equality case of the inequality (\ref{ginoux-1}). It is obvious to check that totally umbilical round spheres in the hyperbolic space $\mathbb{H}^{n+1}$ satisfy the equality in this estimate, however, it is still unknown if they are the only ones. In fact, if the hypersurface is embedded, this result easily follows from the Alexandrov theorem in the hyperbolic space (see \cite{Montiel99}). However, if the hypersurface is only assumed to be {\it immersed} the question is still open. In order to settle this problem, we adopt a method introduced by O.~Hijazi and S.~Montiel in \cite{HijaziMontielhypsphere} which relies on the fact that such hypersurfaces are \emph{critical points} for some eigenvalue functional associated to some Dirac-type operator on $M$. The main result of this section concerns the case when $\wit{M}=\mathbb{H}^{n+1}$ but actually we will prove the following more general statement:

\begin{ethm}\label{t:eqcasegeneral}
Let $M^n$ be an oriented, compact and connected hypersurface immersed into a Riemannian spin manifold $(\wit{M}^{n+1},g)$. If $\wit{M}$ carries a $(i\varepsilon/2)$-Killing spinor for some $\varepsilon\in\{\pm1\}$, then {\rm (\ref{ginoux-1})} (as well as {\rm (\ref{eq:majlambda1L})}) holds and if equality holds then the mean curvature $H$ is constant. Moreover, if $\wit{M}$ also carries a $(-i\varepsilon/2)$-Killing spinor, then equality holds if and only if $M$ is totally umbilical with constant mean curvature.
\end{ethm}

Since the standard hyperbolic space $\mathbb{H}^{n+1}$ has both $(i/2)$- {\it and} $(-i/2)$-Killing spinors (see e.g. \cite{Baum89b}), the previous result immediately implies 
\begin{ecor}\label{t:eqcasehypspace}
The only oriented, compact and connected hypersurfaces immersed into the hyperbolic space $\mathbb{H}^{n+1}$ satisfying $\lambda_1(D_M^2)=(n^2/4)(H^2-1)$ are the totally umbilical round spheres.
\end{ecor}

In Section \ref{PseudoHypSpace}, we will discuss the case of pseudo-hyperbolic spaces.


\subsection{The Hijazi-Montiel approach in presence of imaginary Killing spinors}


Assume that the ambient manifold $\wit{M}$ carries a $(i/2)$-Killing spinor $\Psi\in\Gamma(\Sigma\wit{M})$. After restriction to $M$, it is a straightforward computation to show that $\Psi$ satisfies the modified Dirac equation
\begin{eqnarray}\label{eq:DiracEqMod}
D_+\Psi=\frac{n}{2}H\Psi
\end{eqnarray} 
where $D_+$ is a zero order modification of the extrinsic Dirac operator defined by
\begin{eqnarray}\label{d:DiracMod}
D_+\varphi:=D\varphi-\frac{n}{2} i\nu\cdot\varphi
\end{eqnarray}
for $\varphi\in\Gamma(\Sigma)$. Note that we do not assume that the mean curvature $H$ is constant for the moment. Suppose however that $H$ is positive everywhere on $M$ and consider the metric conformally related to $g$ on $M$ defined by $\ovl{g}:=H^2g$. It is a well-known fact (see \cite{Hitchin,HijaziConf}) that under a conformal change of the metric, there exists a bundle isometry $\varphi\mapsto\ovl{\varphi}$, $\Sigma\to\ovl{\Sigma}$, between the two extrinsic spinor bundles $\Sigma$ and $\ovl{\Sigma}$ over $(M^n,g)$ and $(M^n,\ovl{g})$. Under this identification, the extrinsic Dirac operators $D$ and $D^H$ associated to $g$ and $\ovl{g}$ and acting respectively on $\Sigma$ and $\ovl{\Sigma}$ are related by 
\begin{eqnarray}\label{eq:DiracConf}
D^H\ovl{\varphi}=H^{-\frac{n+1}{2}}\ovl{D(H^{\frac{n-1}{2}}\varphi)}
\end{eqnarray}   
for all $\varphi\in\Gamma(\Sigma)$. Now consider on $\ovl{\Sigma}$ the zero order modification of the extrinsic Dirac operator $D^H$ given by
\begin{eqnarray*}
D^H_+\ovl{\varphi}:=D^H\ovl{\varphi}-\frac{n}{2} H^{-1}\mathcal{I}_\nu\ovl{\varphi}
\end{eqnarray*}
where $\mathcal{I}_\nu$ is the Hermitian endomorphism of $\ovl{\Sigma}$ defined by $\mathcal{I}_\nu\ovl{\varphi}:=\ovl{i\nu\cdot\varphi}$ for all $\varphi\in\Gamma(\Sigma)$. Notice that $D^H_+$ is an elliptic and self-adjoint differential operator of order one which, since $M$ is assumed to be compact, has a discrete spectrum. In the following, we will denote by $\lambda_1(D^H_+)$ the first non-negative eigenvalue of $D^H_+$. Now for every $\varphi\in\Gamma(\Sigma)$, consider the spinor field $\varphi_H:=H^{-\frac{n-1}{2}}\varphi\in\Gamma(\Sigma)$ which is easily seen to satisfy
\begin{eqnarray*}
D^H_+\ovl{\varphi}_H=H^{-\frac{n+1}{2}}\ovl{D_+\varphi}
\end{eqnarray*}
using the conformal covariance (\ref{eq:DiracConf}) of $D$. Taking the $(i/2)$-Killing spinor $\Psi\in\Gamma(\Sigma\wit{M})$ in the previous identity and using (\ref{eq:DiracEqMod}) give that $D^H_+\ovl{\Psi}_H=\frac{n}{2}\ovl{\Psi}_H$. This immediately implies that $\lambda_1(D^H_+)\leq\frac{n}{2}$. Furthermore, if the mean curvature $H$ is constant, it is an easy computation using $\{D,i\nu\cdot\}=0$ to show that 
$${\rm Spec}\big((D^H_+)^2\big)=\Big\{\lambda_{k}\big((D^H_+)^2\big)=H^{-2}\big(\lambda_k(D)^2+(n^2/4)\big)\,/\,\lambda_k(D)\in{\rm Spec}(D)\Big\},$$
so that $\lambda_1(D^H_+)=\frac{n}{2}$ if and only if $\lambda_1(D^2)=\frac{n^2}{4}\left(H^2-1\right)$. Thus we have proved

\begin{prop}\label{CritCar}
Let $M$ be an orientable, compact and connected hypersurface immersed in a Riemannian spin manifold
$(\wit{M}^{n+1},g)$ admitting a $(i/2)$-Killing spinor and suppose that the mean curvature of $M$, after a suitable choice of the unit normal, satisfies $H>0$. Then the first non-negative eigenvalue of $D^{H}_+$ satisfies $\lambda_1(D^H_+)\leq\frac{n}{2}$. Moreover, if $H$ is constant, equality occurs if and only if equality occurs in {\rm (\ref{ginoux-1})}.
\end{prop}

From this proposition, we deduce that any immersion for which (\ref{ginoux-1}) (or equivalently (\ref{eq:majlambda1L})) is an equality realizes a \emph{maximum} for the map 
\begin{eqnarray*}
\mathcal{F}^+_1:\iota\in{\rm Imm}^+(M,\wit{M})\mapsto\lambda_1(D^{H_\iota}_+)\in\mathbb{R}
\end{eqnarray*}
where ${\rm Imm}^+(M,\wit{M})$ denotes the space of isometric immersions of $M$ in $\wit{M}$ with non-vanishing mean curvature $H_\iota$. This characterization of hypersurfaces satisfying the equality case in (\ref{ginoux-1}) leads to the study of the critical points of the functional $\mathcal{F}^+_1$. 

\begin{erem}
{\rm 
It is important to note that if the manifold $\wit{M}$ carries a $(-i/2)$-Killing spinor, then Proposition \ref{CritCar} is true with the operators $D_+$ and $D_+^H$ replaced respectively by 
\begin{eqnarray*}
D_-:=D+\frac{n}{2}i\nu\cdot:\Gamma(\Sigma)\rightarrow\Gamma(\Sigma)
\end{eqnarray*}
and
\begin{eqnarray}
D^H_-:=D^H+\frac{n}{2}H^{-1}\mathcal{I}_\nu:\Gamma(\ovl{\Sigma})\rightarrow\Gamma(\ovl{\Sigma}).
\end{eqnarray}
In this situation, the corresponding functional is defined by
\begin{eqnarray*} 
\mathcal{F}^-_1:\iota\mapsto\lambda_1^-(D^{H_\iota}_-) 
\end{eqnarray*}
where $\lambda_1^-(D^{H_\iota}_-)$ is the first non-negative eigenvalue of $D^{H_\iota}_-$.}
\end{erem}


\subsection{Derivatives of the functional $\mathcal{F}^\pm_1$ }


As explained in the previous section we are led to study the first derivatives of the functional $\mathcal{F}^\pm_1$ at least in a particular situation. As above, we start with an immersion $\iota=\iota_0:M\to\wit{M}$ with positive mean curvature (not necessarily constant) and such that $\lambda_1(D^H_+)=\frac{n}{2}$. Note that here {\it we do not assume the existence of imaginary Killing spinor fields on $\wit{M}$}.

Now we deform the immersion $\iota$ along normal geodesics, that is, we consider, for $\varepsilon>0$ sufficiently small, the map $F:]-\varepsilon,\varepsilon[\times M\to\wit{M}$, $(t,x)\mapsto\exp_{\iota(x)}(t\nu_x)$. Note that, choosing $\varepsilon>0$ sufficiently small, the map $F$ is smooth and $F(t,\cdot):M\to\wit{M}$ is an immersion such that $F(0,\cdot)=\iota$. In fact, the map $t\mapsto F(t,x)$ is the geodesic starting from $\iota(x)$ with speed vector $\nu_x$, and so it is analytic. For each $t\in]-\varepsilon,\varepsilon[$, we denote by $g_t:=F(t,\cdot)^*g$ the induced metric on $M$, by $\nu_t$ the unit normal field inducing the orientation of $M$, by $H_t:=-(1/n){\rm tr}(\wnabla\nu_t)$ the mean curvature of $F(t,\cdot)$ -- which, up to making $\varepsilon>0$ smaller, may be assumed to be positive on $M$ for all $t\in]-\varepsilon,\varepsilon[$ -- and by $\ovl{g}_t:=H_t^2 g_t$. We also denote by $D^{H_t}$ the Dirac operator associated to the metric $\ovl{g}_t$ and let $D^{H_t}_+:=D^{H_t}-\frac{n}{2}H_t^{-1}\mathcal{I}_{\nu_t}:\Gamma(\ovl{\Sigma}_{t})\to\Gamma(\ovl{\Sigma}_t)$, where $\mathcal{I}_{\nu_t}$ is the Hermitian endomorphism of $\ovl{\Sigma}_t$ defined by $\mathcal{I}_{\nu_t}\ovl{\varphi}:=\ovl{i\nu_t\cdot\varphi}$. Here $\ovl{\Sigma}_t$ denotes the extrinsic spinor bundle over $M$ endowed with the spin structure induced by $\wit{M}$ and the Riemannian metric $\ovl{g}_t$. Since we perturb the immersion analytically, the family $(D^{H_t}_+)$ with ${t\in]-\varepsilon,\varepsilon[}$ is an analytic family of unbounded closed self-adjoint operators with compact resolvent, therefore the spectrum of $D^{H_t}_+$ can be written as a sequence $(\mu^+_k(t))_{k\in\mathbb{N}}$, where each eigenvalue $\mu_k^+(t)$ depends analytically on $t$ and where corresponding eigenvectors can be found to also depend analytically on $t$ (see \cite{Katobook}). We denote by $\lambda_1^+(t)$ any branch of that spectrum with $\lambda_1^+(0)=\lambda_1(D^{H}_+)$, the smallest non-negative eigenvalue of $D^H_+=D^{H_0}_+$. Following \cite{BaerGaudMor}, we denote by $\tau_0^t:\ovl{\Sigma}_0=\ovl{\Sigma}\to\ovl{\Sigma}_t$ the pa\-ral\-lel transport along the curves $s\mapsto(s,x)$ in the so-called generalized cylinder $\left(]-\varepsilon,\varepsilon[\times M,dt^2\oplus\ovl{g}_t\right)$, for all $t\in]-\varepsilon,\varepsilon[$. Then for any analytic family $(\ovl{\Phi}_t)_t$ of eigenvectors associated to $\lambda_1^+(t)$, differentiating the identity 
\[\lambda_1^+(t)\int_M|\ovl{\Phi}_t|^2dv_{\ovl{g}_t}=\int_M\Re e\langle D^{H_t}_+\ovl{\Phi}_t,\ovl{\Phi}_t\rangle dv_{\ovl{g}_t}\]
at $t=0$ yields
\[\frac{d\lambda_1^+}{dt}(0)\int_M|\ovl{\Phi}_0|^2dv_{\ovl{g}_0}=\int_M\Re e\langle\frac{d}{dt}\Big|_{t=0}\pa{\tau_t^0D^{H_t}_+\tau_0^t\ovl{\Phi}_0},\ovl{\Phi}_0\rangle dv_{\ovl{g}_0}.\]
Now we have $\tau_t^0D^{H_t}_+\tau_0^t=\tau_t^0 D^{H_t}\tau_0^t-\frac{n}{2}H_t^{-1}\tau_t^0\mathcal{I}_{\nu_t}\tau_0^t$ and, since the variation of $\iota$ is a geodesic normal one, the vector field $\nu_t=\frac{\partial}{\partial t}$ is parallel along the curves $s\mapsto(s,x)$, so that $\tau_t^0\mathcal{I}_{\nu_t}\tau_0^t=\mathcal{I}_{\nu_0}=\mathcal{I}_{\nu}$ for all $t\in]-\varepsilon,\varepsilon[$.
With the formula for the first variation of the Dirac operator by J.-P.~Bourguignon and P.~Gauduchon \cite{BourgGaud92} (see also \cite{BaerGaudMor}), we deduce that
\be 
\frac{d\lambda_1^+}{dt}(0)\int_M|\ovl{\Phi}_0|^2dv_{\ovl{g}_0} &=&-\frac{1}{2}\int_M\ovl{g}_0\pa{T_{\ovl{\Phi}_0},\frac{\partial\ovl{g}_t}{\partial t}(0)}dv_{\ovl{g}_0}\\
&&+\frac{n}{2}\int_MH^{-2}\frac{\partial H_t}{\partial t}\Big|_{t=0}\Re e\langle\mathcal{I}_\nu\ovl{\Phi}_0,\ovl{\Phi}_0\rangle dv_{\ovl{g}_0},
\ee
where 
\begin{eqnarray*}
T_{\ovl{\Phi}_0}(X,Y):=\frac{1}{2}\Re e\langle X\ovl{\cdotSi}\ovl{\nabla}_Y\ovl{\Phi}_0+Y\ovl{\cdotSi}\ovl{\nabla}_X\ovl{\Phi}_0,\ovl{\Phi}_0\rangle
\end{eqnarray*}
is the so-called \emph{energy-momentum} tensor associated to $\ovl{\Phi}_0$. Here $\ovl{\cdotSi}$ is the Clifford multiplication on $\ovl{\Sigma}$ defined by (\ref{d:CliffMult}) and $\ovl{\nabla}$ is the spin Levi-Civita connection with respect to the metric $\ovl{g}_0$. Note that we kept the same notations for the Hermitian scalar products on $\ovl{\Sigma}$ and $\Sigma$. Now fix an eigenvector $\ovl{\Phi}_0\in\Gamma(\ovl{\Sigma})$ for the Dirac-type operator $D^{H}_+$ associated with $\lambda_1(D^H_+)$ and let $\ovl{\Psi}_0:=H^{\frac{n-1}{2}}\ovl{\Phi}_0$.  We compute $\frac{d\lambda_1^+}{dt}(0)$ in terms of $\Psi_0\in\Gamma(\Sigma)$ and of geometric quantities attached to $\iota$. First, since $\frac{\partial F}{\partial t}(0,\cdot)=\nu$, we have on the one hand (see e.g. \cite{Montiel99})
\be 
\frac{\partial\ovl{g}_t}{\partial t}(0)=\frac{\partial}{\partial t}\Big|_{t=0}\pa{H_t^2g_t}=\frac{2H}{n}\pa{|A|^2+\wit{{\rm ric}}(\nu,\nu)} g-2H^2g(A\cdot,\cdot).
\ee
On the other hand, using the isomorphism $\Sigma\to\ovl{\Sigma}$, we may write (see e.g. \cite[Sec. 1.3]{Ginouxbook})
\be 
T_{\ovl{\Phi}_0}(X,Y)=H^{-n+2}T_{\Psi_0}(X,Y),
\ee
for all $X,Y\in\Gamma(TM)$, where $T_{\Psi_0}$ is the energy-momentum tensor associated to $\Psi_0$ defined by
\begin{eqnarray*}
T_{\Psi_0}(X,Y):=\frac{1}{2}\Re e\langle X\cdotSi\nabla_Y\Psi_0+Y\cdotSi\nabla_X\Psi_0,\Psi_0\rangle.
\end{eqnarray*}
Therefore, assuming without loss of generalities that $\int_M|\ovl{\Phi}_0|^2dv_{\ovl{g}_0}=1$, we compute:
\begin{eqnarray*}
\frac{d\lambda_1^+}{dt}(0)&=&\frac{1}{n}\int_M H^{-1}\pa{|A|^2+\wit{{\rm ric}}(\nu,\nu)}\left(\frac{n}{2}\Re e\langle i\nu\cdot\Psi_0,\Psi_0\rangle-g(T_{\Psi_0},g)\right) dv_{g}\\
& & +\int_Mg(T_{\Psi_0},A)dv_{g}.
\end{eqnarray*}
But since $g(T_{\Psi_0},g)={\rm tr}_g(T_{\Psi_0})=\Re e\langle D\Psi_0,\Psi_0\rangle$, we obtain
\begin{eqnarray*}
\frac{d\lambda_1^+}{dt}(0) =-\frac{1}{n}\int_M H^{-1}\pa{|A|^2+\wit{{\rm ric}}(\nu,\nu)}\Re e\langle D_+\Psi_0,\Psi_0\rangle dv_{g}+\int_Mg(T_{\Psi_0},A)dv_{g}.
\end{eqnarray*}
However, since $\ovl{\Phi}_0\in\Gamma(\ovl{\Sigma})$ is an eigenspinor for $D^H_+$ associated with the eigenvalue $\lambda_1^+(0)=\frac{n}{2}$ and from the equivalence
\begin{eqnarray}\label{EquivConfEigen}
D^H_+\ovl{\Phi}_0=\frac{n}{2}\ovl{\Phi}_0\quad\Longleftrightarrow\quad D_+\Psi_0=\frac{n}{2}H\Psi_0,
\end{eqnarray}
one concludes that
\begin{eqnarray}\label{FirstStep}
\frac{d\lambda_1^+}{dt}(0) =-\frac{1}{2}\int_M \pa{|A|^2+\wit{{\rm ric}}(\nu,\nu)}|\Psi_0|^2 dv_{g}+\int_Mg(T_{\Psi_0},A)dv_{g}.
\end{eqnarray}
To compute the remaining term $g(T_{\Psi_0},A)$, we define a new covariant derivative by $\wih{\nabla}^+_X:=\wnabla_X-(i/2)X\cdot$ on $\Sigma$. Then a lengthy but direct calculation using the spin Gau\ss{} formula (\ref{eq:Gaussspin}) yields that for any $\varphi\in\Gamma(\Sigma)$,
\be 
|\wih{\nabla}^+\varphi|^2&:=&\sum_{j=1}^n|\wih{\nabla}^+_{e_j}\varphi|^2\\
& = &\sum_{j=1}^n|\nabla_{e_j}\varphi+\frac{A(e_j)}{2}\cdot\nu\cdot\varphi-\frac{i}{2}e_j\cdot\varphi|^2\\
&=&|\nabla\varphi|^2+\left(\frac{|A|^2+n}{4}\right)|\varphi|^2-g(T_{\varphi},A)-\Re e\langle i\nu\cdot (D\varphi-\frac{nH}{2}\varphi),\varphi\rangle.
\ee
For $\varphi=\Psi_0$, we deduce using the right-hand side of (\ref{EquivConfEigen}) that 
\be 
g(T_{\Psi_0},A)=|\nabla\Psi_0|^2-|\wih{\nabla}^+\Psi_0|^2+\left(\frac{|A|^2-n}{4}\right)|\Psi_0|^2.
\ee
Now integrating over $M$ this identity with the help of the famous Schr\"odinger-Lichnerowicz-formula
\be
D^2=\nabla^*\nabla+\frac{S}{4}
\ee
gives 
\be 
\int_Mg(T_{\Psi_0},A)dv_{g}=\int_M\pa{\Re e\langle D^2\Psi_0,\Psi_0\rangle-\frac{S}{4}|\Psi_0|^2-|\wih{\nabla}^+\Psi_0|^2+\left(\frac{|A|^2-n}{4}\right)|\Psi_0|^2}dv_{g}.
\ee
Here $S$ stands for the scalar curvature of $(M^n,g)$. On the other hand, from (\ref{eq:commutDf}), (\ref{EquivConfEigen}) and the anti-commutativity rule $\{D,i\nu\cdot\}=0^+$, we check that
\be 
D^2\Psi_0=\frac{n^2}{4}(H^2-1)\Psi_0+\frac{n}{2}\nabla H\cdot\nu\cdot\Psi_0,
\ee
so that $\Re e\langle D^2\Psi_0,\Psi_0\rangle=\frac{n^2}{4}(H^2-1)|\Psi_0|^2$ and hence
\[\int_Mg(T_{\Psi_0},A)dv_{g}=\int_M\pa{\frac{1}{4}\left(n^2(H^2-1)-S+|A|^2-n\right)|\Psi_0|^2-|\wih{\nabla}^+\Psi_0|^2}dv_g.\]
The Gau\ss{} formula for the scalar curvature provides 
\[S=\wit{S}-2\wit{{\rm ric}}(\nu,\nu)+n^2H^2-|A|^2,\]
from which
\be 
\int_Mg(T_{\Psi_0},A)dv_{g} & = & -\int_M\left(\frac{1}{4}\left(\wit{S}+n(n+1)\right)-\frac{1}{2}\left(|A|^2+\wit{{\rm ric}}(\nu,\nu)\right)\right)|\Psi_0|^2dv_g\\ 
& & -\int_M|\wih{\nabla}^+\Psi_0|^2dv_g
\ee
follows. Inserting this identity in (\ref{FirstStep}), we finally deduce that 
\be 
\frac{d\lambda^+_1}{dt}(0)=-\int_M\pa{|\wih{\nabla}^+\Psi_0|^2+\frac{\wit{S}+n(n+1)}{4}|\Psi_0|^2}dv_g.
\ee
It is worth noticing that this formula holds if we assume that it is the first non-negative eigenvalue $\lambda_1(D^H_-)$ of $D^H_-$ which satisfies $\lambda_1(D^H_-)=\frac{n}{2}$ instead of $\lambda_1(D^H_+)$; in this situation, $\wih{\nabla}^+$ has to be replaced with the covariant derivative  defined by $\wih{\nabla}^-_X:=\wit{\nabla}_X+(i/2)X\cdot$. 

From this computation, it is now straightforward to give a necessary condition for an immersion $\iota$ to be a critical point of $\mathcal{F}_1^\pm$:
\begin{ethm}\label{t:critical}
Let $M$ be an oriented, compact and connected hypersurface isometrically immersed in a Riemannian spin manifold $(\wit{M}^{n+1},g)$. Assume that the scalar curvature $\wit{S}$ of $\wit{M}$ is greater or equal to $-n(n+1)$ and that the mean curvature $H$ of $M$ with respect to a suitable choice of the normal is positive. If $\lambda_1(D^H_\varepsilon)=\frac{n}{2}$ for some $\varepsilon\in\{\pm1\}$ and it is critical for all the variations of the hypersurface $M$ in $\wit{M}$, then $\wit{S}=-n(n+1)$ and $\wit{\nabla}_X\Psi=(i\varepsilon/2)X\cdot\Psi$ for all $X\in\Gamma(TM)$ for all $\Psi\in\Gamma(\Sigma)$ satisfying
\begin{eqnarray*}
D_\varepsilon\Psi=\frac{n}{2}H\Psi.
\end{eqnarray*} 
\end{ethm}


\subsection{Proof of Theorem \ref{t:eqcasegeneral}}


If $\wit{M}$ carries a $(i\varepsilon/2)$-Killing spinor for some $\varepsilon\in\{\pm1\}$, then from Theorem \ref{t:lambda1L} and Remark \ref{r:knownestimates} the inequalities (\ref{eq:majlambda1L}) and (\ref{ginoux-1}) hold. Moreover, if equality holds in (\ref{eq:majlambda1L}), Proposition \ref{p:equality} implies that the mean curvature is constant and then $\lambda_1(D)^2=\frac{n^2}{4}(H^2-1)$.

Assume now that $\wit{M}$ carries a $(i/2)$- as well as a $(-i/2)$-Killing spinor. From Proposition \ref{CritCar}, we deduce that such an immersion is a maximum for the functional $\mathcal{F}^+_1$ and thus $\frac{d\lambda_1^+}{dt}(0)=0$. Let $\Phi$ be a non-zero $(-i/2)$-Killing spinor on $\wit{M}$ so that $D_-\Phi=\frac{n}{2}H\Phi$. From this equation and since $H$ is constant, a direct computation shows that the spinor $\wit{\Phi}:=H\Phi-i\nu\cdot\Phi$ satisfies $D_+\wit{\Phi}=\frac{n}{2}H\wit{\Phi}$.  
On the other hand, since the existence of an $(\pm i/2)$-Killing spinor on $\wit{M}$ implies that $\wit{M}$ is an Einstein manifold with scalar curvature $\wit{S}=-n(n+1)$ (see \cite{BFGK} for example), Theorem \ref{t:critical} applies and we get that $\wnabla_X\wit{\Phi}=(i/2)X\cdot\wit{\Phi}$ for all $X\in\Gamma(TM)$, that is 
\be 
\frac{i}{2}X\cdot\pa{H\Phi-i\nu\cdot\Phi}&=&\wnabla_X\pa{H\Phi-i\nu\cdot\Phi}\\
&=&H\pa{-\frac{i}{2}X\cdot\Phi}+iA(X)\cdot\Phi-i\nu\cdot\pa{-\frac{i}{2}X\cdot\Phi}\\
&=&iA(X)\cdot\Phi-\frac{iH}{2}X\cdot\Phi-\frac{i}{2}X\cdot i\nu\cdot\Phi.
\ee
This implies that $(A(X)-HX)\cdot\Phi=0$ for all $X\in\Gamma(TM)$, and since $\Phi$ has no zero, $M$ is totally umbilical. This concludes the proof of Theorem \ref{t:eqcasegeneral}.


\subsection{The case of pseudo-hyperbolic spaces}\label{PseudoHypSpace}


In this section, we examine the case of other {\it complete} ambient manifolds $\wit{M}$ carrying imaginary Killing spinors. These manifolds have been classified by H.~Baum \cite{Baum89a,Baum89b} and are known as {\it pseudo-hyperbolic spaces}.
For the sake of completeness and since we need an additional argument for our purpose, we recall the result of \cite{Baum89a,Baum89b} and give a sketch of the proof:

\begin{prop}\label{p:classifBaum}
Let $(\wit{M}^{n+1},g)$ be a complete Riemannian spin manifold admitting a non-zero $(i\varepsilon/2)$-Killing spinor for some $\varepsilon\in\{\pm1\}$. Then $(\wit{M}^{n+1},g)$ is isometric to either the real hyperbolic space of constant sectional curvature $-1$ or to the warped product $(\R\times N,dt^2\oplus e^{2 t}g_N)$, where $(N^n,g_N)$ is a complete non-flat Riemannian spin manifold carrying at least one non-zero parallel spinor. In the latter case, denoting by $\mathcal{K}_0(N,g_N)$ (resp. $\mathcal{K}_0^\varepsilon(N,g_N)$) the space of parallel spinors on $(N^n,g_N)$ for the induced metric and spin structure (resp. its projection onto the half-spinors bundle $\Sigma_\varepsilon N$ if $n$ is even), the map
\be
\left.
\begin{array}{ll}
\mathcal{K}_0^\varepsilon(N,g_N) & \textrm{ if }n\textrm{ is even}\\ 
\mathcal{K}_0(N,g_N) & \textrm{ if }n\textrm{ is odd}
\end{array}
\right| & \lra&\left\{\frac{i\varepsilon}{2}-\textrm{Killing spinors on }\wit{M}\right\}\\
\varphi & \lmt & \left|
\begin{array}{ll}
e^{\frac{t}{2}}\varphi&\textrm{ if }n\textrm{ is even}\\
e^{\frac{t}{2}}(\varphi\oplus\varepsilon i\frac{\partial}{\partial t}\cdot\varphi) & \textrm{ if }n\textrm{ is odd,}
\end{array}
\right.
\ee
is a well-defined monomorphism. If moreover $N$ is \emph{compact}, then this is actually an isomorphism.
\end{prop}

{\it Proof}: Let $\varphi$ be a non-zero $(i\varepsilon/2)$-Killing spinor on the manifold $(\wit{M}^{n+1},g)$. As H.~Baum showed (see \cite{Baum89a} and references therein), if $(\wit{M},g)$ is \emph{not} isometric to the hyperbolic space, then there must exist an unit smooth vector field $\xi$ on $\wit{M}$ with $i\xi\cdot\varphi=\varepsilon\varphi$ on $\wit{M}$. From this relationship, the foliated structure of $\wit{M}$ can be elementary deduced as follows. First note that $\xi=(\varepsilon V)/|V|$, where $g(V,X):=i\langle X\cdot\varphi,\varphi\rangle$ for all $X\in\Gamma(T\wit{M})$ and in particular $V=\varepsilon\nabla|\varphi|^2$ has no zeros on $\wit{M}$. Since $\wnabla_XV=\varepsilon|\varphi|^2X$ (that is $V$ is a closed conformal vector field on $\wit{M}$), one deduces that $\wnabla_X\xi=X-g(X,\xi)\xi$ for all $X\in\Gamma(T\wit{M})$ and as a consequence, the flow of $\xi$, which is well-defined and complete since $(\wit{M},g)$ is complete, preserves the level hypersurfaces of $|\varphi|^2=|V|$. On the other hand, the second fundamental form of each such hypersurface with respect to $\xi$ is $-\mathrm{Id}$, the Lie derivative of the metric in the direction of $\xi$ is given by $\mathcal{L}_\xi g=2g_{|_{\xi^\perp\times\xi^\perp}}$ and hence, setting 
\begin{eqnarray*}
N:=\left\{x\in\wit{M}\,,\,|\varphi|^2(x)=1\right\}\subset\wit{M}, 
\end{eqnarray*}
the flow of $\xi$ provides a diffeomorphism $\R\times N\to\wit{M}$ identifying $\xi$ with $\frac{\partial}{\partial t}$ and pulling back the metric $g$ onto $dt^2\oplus e^{2t}g_N$, where $g_N$ is the metric induced from $g$ onto $N$. This done, the spin Gau\ss{} formula (\ref{eq:Gaussspin}) implies that, for any $X\in\Gamma(TN)$,
\[\frac{i\varepsilon}{2}X\cdot\varphi=\wnabla_X\varphi=\nabla_X^{\Sigma N}\varphi-\frac{X}{2}\cdot\xi\cdot\varphi=\nabla_X^{\Sigma N}\varphi+\frac{i\varepsilon}{2}X\cdot\varphi,\]
from which $\nabla^{\Sigma N}\varphi_{|_N}=0$ follows: the restriction of $\varphi$ onto any level hypersurface of $|\varphi|^2$ is a parallel spinor. Here $\nabla^{\Sigma N}$ stands for the spin Levi-Civita connection on $\Sigma:=\Sigma\widetilde{M}_{|N}$. In case $n$ is even, the condition $i\xi\cdot\varphi=\varepsilon\varphi$ actually imposes $\varphi\in\Gamma(\Sigma_\varepsilon N)$ since $i\xi\cdot$ coincides with the Clifford action of the complex volume form of $(N,g_N)$.
In case $n$ is odd, the spinor $\varphi_{|_N}$ can be rewritten in the form $\varphi_{|_N}=\varphi_0\oplus\varepsilon i\frac{\partial}{\partial t}\cdot\varphi_0$, where $\varphi_0\in\Gamma(\Sigma N)$ is parallel.
The dependence in $t$ of $\varphi$ is easily computed thanks to
\[\frac{\partial\varphi}{\partial t}=\wnabla_{\frac{\partial}{\partial t}}\varphi=\frac{i\varepsilon}{2}\frac{\partial}{\partial t}\cdot\varphi=\frac{1}{2}\varphi,\]
from which $\varphi(t,\cdot)=e^{\frac{t}{2}}\varphi(0,\cdot)$ follows.
This gives the formulas for the above map, which is obviously a right inverse to the ``restriction'' map
\be \left\{\frac{i\varepsilon}{2}-\textrm{Killing spinors on }\wit{M}\right\}&\lra&\left|\begin{array}{ll}\mathcal{K}_0^\varepsilon(N,g_N)&\textrm{ if }n\textrm{ is even}\\\mathcal{K}_0(N,g_N)&\textrm{ if }n\textrm{ is odd}\end{array}\right.\\
\varphi&\lmt&\left|\begin{array}{ll}\varphi_{|_{\{0\}\times N}}&\textrm{ if }n\textrm{ is even}\\\varphi_+{}{|_{\{0\}\times N}}&\textrm{ if }n\textrm{ is odd.}\end{array}\right.
\ee
In case $N$ is compact, this restriction map is surjective, a remark missing in \cite{Baum89b}. To show this, let $\psi$ be any further non-zero $(i\varepsilon/2)$-Killing spinor on $(\wit{M}^{n+1},g)$.
Then, again, $\psi$ splits $(\wit{M}^{n+1},g)$ as a warped product $(\R\times P^n,ds^2\oplus e^{2s}g_\Sigma)$, where $(P^n,g_P)$ is complete, spin and carries a non-zero parallel spinor. 
Now, using  the work \cite{Montiel99} of S.~Montiel, the latter splitting must ``coincide'' (in a sense that is made precise below) with the former. Namely, for all $t\in\R$ the hypersurface $\{t\}\times N$ is a totally umbilical compact hypersurface of $\wit{M}$ with constant mean curvature. Therefore, by applying \cite[Lemma 4]{Montiel99} to the foliation of $\wit{M}$ induced by $\psi$ (whose leaves are not assumed to be compact), we easily conclude that for each $t\in\R$, there exists an $s\in\R$ such that $\{t\}\times N=\{s\}\times P$; in particular, $P$ itself must be compact. The same argument shows that, for each $s\in\R$, there exists a $t\in\R$ with $\{s\}\times P=\{t\}\times N$. This yields that, if $\Phi:\R\times P\lra\R\times N$, $(s,x)\mapsto(\phi_1(s,x),\phi_N(s,x))$, is the isometry induced by both splittings, then the component map $\phi_1$ already only depends on $s$. By $\Phi^*(dt^2\oplus e^{2t}g_N)=ds^2\oplus e^{2s}g_P$ and the existence of an inverse map for $\Phi$ of a similar form, one deduces on the one hand that $\frac{\partial\phi_N}{\partial s}(s,x)=0$ and hence $(\phi_1'(s))^2=1$ for all $s\in\R$, and on the other hand that $e^{2s}g_P=e^{2\phi_1(s)}(\phi_N)^*g_N$ holds for all $s\in\R$. This in turn implies the existence of an $s_0\in\R$ with $\phi_1(s)=s-s_0$ and $g_P=e^{-2s_0}(\phi_N)^*g_N$. Thus, up to homotheties on the metrics $g_P$ and $g_N$, the Riemannian manifolds $(P,g_P)$ and $(N,g_N)$ are isometric and, up to translations in $s$, the splittings $\R\times P$ and $\R\times N$ coincide.
By the first part of the proof, $\psi$ must come from a parallel spinor on $N$ and hence lie in the image of the map of Proposition \ref{p:classifBaum}. This concludes the proof.
\findemo

From the previous result, we deduce a characterization of hypersurfaces for which Inequality (\ref{eq:majlambda1L}) is an equality when $\wit{M}$ is a pseudo-hyperbolic space in several situations. In fact, as we will see, we are left with the case $n$ is even, the manifold $(N^n,g_N)$ has only positive (or only negative) non-zero parallel spinors and $M$ is only immersed in $\wit{M}$. Indeed, we prove
\begin{ecor}\label{c:Baumclass}
Let $(\wit{M}^{n+1},g):=(\R\times N,dt^2\oplus e^{2t}g_N)$, where $(N^n,g_N)$ is a closed non-flat Riemannian spin manifold endowed with at least one non-zero pa\-ral\-lel spinor and assume that $\wit{M}$ carries the induced spin structure (in particular, $(\wit{M},g)$ admits an imaginary Killing spinor for at least one of the constants $(\pm i/2)$). Let $M^n\hookrightarrow\wit{M}$ be any immersed closed orientable hypersurface carrying the induced metric and spin structure and suppose that one of the following supplementary conditions is fulfilled:
\beit\item[a)] $n$ is odd;
\item[b)] $n$ is even and $(N^n,g_N)$ has non-zero positive as well as negative pa\-ral\-lel spinors;
\item[c)] $n$ is even and $M^n$ bounds a domain in $\wit{M}$.
\eeit
Then $M$ satisfies the equality case in {\rm(\ref{eq:majlambda1L})} (and so in {\rm(\ref{ginoux-1})}) if and only if $M=\{t\}\times N$ for some $t\in\R$.
\end{ecor}

{\it Proof}:  
From Proposition \ref{p:equality}, we have that if $M^n\hookrightarrow\wit{M}^{n+1}$ satisfies the equality case in (\ref{eq:majlambda1L}) then its mean curvature $H$ must be constant. If either $a)$ or $b)$ is fulfilled, then by Proposition \ref{p:classifBaum}, the manifold $(\wit{M},g)$ admits non-zero imaginary Killing spinors for \emph{both} constants $(\pm i/2)$, therefore Theorem \ref{p:equality} implies that $M$ is totally umbilical which, combined with \cite[Lemma 4]{Montiel99}, yields $M=\{t\}\times N$ for some $t\in\R$. If $c)$ is fulfilled, this time \cite[Theorem 10]{Montiel99} applies and yields again $M=\{t\}\times N$ for some $t\in\R$. This shows the ``only if'' part of the corollary. The ``if'' part is easy to see since $\lambda_1(D_M)=0$ because of parallel spinors on $N$, and on the other hand $|H|=1$ by the explicit form of the metric. This concludes the proof.
\findemo

\vspace{0.5cm}

{\bf Acknowledgment.} Part of this work was carried out during the Research in Pairs no. 998 at the Centre International de Rencontres Math\'ematiques (Luminy, France), which the first two named authors would like to thank for their generous support and friendly welcome.
Discussions at a latter stage of the paper were made possible thanks to a travel grant awarded to the first named author by the German Academic Exchange Service (DAAD) and the hospitality of the Lebanese University at Beirut, that the first named author would like to thank.

\providecommand{\bysame}{\leavevmode\hbox to3em{\hrulefill}\thinspace}

\end{document}